\documentclass[12pt]{amsart}

\newtheorem{thm}{Theorem}[section]   % Numbered within each section
\newtheorem{cor}[thm]{Corollary}     % Numbered along with thm
\newtheorem{lem}[thm]{Lemma}         % Numbered along with thm
\newtheorem{prop}[thm]{Proposition}  % Numbered along with thm

\theoremstyle{definition}
\newtheorem{defn}[thm]{Definition}   % Numbered along with thm

\theoremstyle{remark}
\newtheorem{ex}[thm]{Example}        % Numbered along with thm

\numberwithin{equation}{section}     % Number equations within sections

\newcommand{\secref}[1]{Section~\ref{#1}}
\newcommand{\thmref}[1]{Theorem~\ref{#1}}
\newcommand{\corref}[1]{Corollary~\ref{#1}}
\newcommand{\lemref}[1]{Lemma~\ref{#1}}
\newcommand{\propref}[1]{Proposition~\ref{#1}}

\renewcommand{\H}{\mathcal H}
\newcommand{\B}{\mathcal B}
\newcommand{\K}{\mathcal K}

\newcommand{\C}{\mathbb C}
\newcommand{\R}{\mathbb R}
\newcommand{\Z}{\mathbb Z}
\newcommand{\N}{\mathbb N}

\newcommand{\T}{\mathbb T}
\newcommand{\cstar}{\text{$C^*$}-}

\newcommand{\spacetext}[1]{\quad\text{#1}\quad}

\newcommand{\noqed}{\renewcommand{\qed}{}}
\newcommand{\ad}{\operatorname{Ad}}
\newcommand{\paut}{\operatorname{PAut}}

\newcommand{\inv}{^{-1}}

\begin{document}

\title{Crossed product duality for partial \cstar automorphisms}
\author{John Quigg}
\address{Department of Mathematics\\Arizona State University\\
 Tempe, Arizona 85287}
\email{quigg@math.la.asu.edu}
\thanks{This material is based upon work supported by the National
Science Foundation under Grant No. DMS9401253}
\date{April 1, 1996 (revised)}

\begin{abstract}
For partial automorphisms of \cstar algebras, Tak\-ai-Tak\-e\-saki crossed 
product duality tends to fail, in proportion to the extent to which the 
partial automorphism is not an automorphism. 
\end{abstract}

\maketitle

\section{Introduction}

Recently Exel \cite{exe:partial} introduced the notion of a crossed 
product of a \cstar algebra by a partial automorphism (an isomorphism 
between ideals), in order to better understand circle actions. This 
generalizes crossed products by automorphisms (equivalently, integer 
actions), and some of the usual theory of crossed products by actions 
carries over to this new context \cite{exe:partial}, \cite{exe:twist}, 
\cite{mcc}, \cite{qui-rae:partial}. It seems natural to ask about the 
Takai-Takesaki crossed product duality \cite{tak}. In this paper we show 
that, perhaps unsurprisingly (since partial automorphisms, being 
partially defined, miss some of the information of the \cstar algebra), 
crossed product duality tends to fail for partial automorphisms. Indeed, 
crossed product duality seems to fail more miserably the more 
``partial'' the partial automorphism is. 

To be more precise, from experience with Takai-Takesaki duality for 
crossed products by abelian groups, we expect a dual action of the 
circle group $\T$ on a crossed product by a partial action, and indeed 
Exel \cite{exe:partial} constructs such a thing. We apologize, but for 
our purposes we find it more convenient to work with the corresponding 
coaction of the integer group $\Z$. For abelian locally compact groups, 
statements about coactions are just Fourier transforms of statements 
about actions of the dual groups. However, at a certain point we need a 
representation of the circle group
$\T$, which is more easily dealt with as a 
representation of $c_0(\Z)$. For the reader's convenience, in 
\secref{setup}, after reviewing the elementary theory of partial 
automorphisms, we give a rough guide to what we need from coactions, 
specialized to $\Z$; suitable references are, e.g., \cite{lprs} and 
\cite{qui:disc}. 

In \secref{invariantideals} we prove what are surely expected results 
concerning invariant ideals, and concerning tensor products with the 
identity automorphism. 

In \secref{duality} we obtain a kind of ``Wold decomposition'' of a 
partial automorphism, showing that there is a largest subalgebra, which 
turns out to be an ideal, on which we have an actual automorphism, and 
the quotient partial automorphism is as far as possible from an 
automorphism (``completely nonautomorphic''). For these completely 
nonautomorphic partial automorphisms, crossed product duality
fails most dramatically, at least under a mild condition that certain
projections be multipliers.

In \secref{distribution} we study how the behavior of the partial action 
$\alpha$ depends upon the distribution of the domains and ranges of the 
powers $\alpha^n$. Along the way see a more fundamental reason why 
crossed product duality
fails in general for partial automorphisms. 

\section{partial automorphisms and the dual coaction}
\label{setup}

We begin with a review and embellishment of some of Exel's work 
\cite{exe:partial}. A {\em partial automorphism\/} of a \cstar algebra 
is an isomorphism between two (closed, two-sided) ideals. Let $\paut A$ 
denote the set of all partial automorphisms of $A$. Since an ideal of an 
ideal of $A$ is an ideal of $A$, $\paut A$ is closed under composition 
(where the domain of $\alpha \circ \beta$ is taken to be all elements in 
the domain of $\beta$ which $\beta$ maps into the domain of $\alpha$). 
For $\alpha \in \paut A$ we let $\alpha^0$ be the identity automorphism 
of $A$, and for $n>0$ we let $\alpha^{-n}$ be the $n$th power of the 
inverse $\alpha\inv$. We denote the range of the partial automorphism 
$\alpha^n$ by $D_n$ (so that the domain is $D_{-n}$). We have 
\[ \dotsb D_{-2} \subset D_{-1} 
  \subset A = D_0 \supset D_1 \supset D_2 \dotsb \]
and
\begin{equation}
\label{mapideals}
\alpha^n(D_kD_{-n}) = D_{n+k}D_n \spacetext{for all} n,k \in \Z.
\end{equation}

We refer to the $D_n$ as the ideals of $\alpha$, and write
$D_n(\alpha)$ when there are more than one partial automorphism
around.  Note that $D_n(\alpha\inv) = D_{-n}(\alpha)$.  Thus, many
properties of the $D_n$ for $n<0$ follow from the corresponding
properties for $n>0$ by replacing $\alpha$ by $\alpha\inv$; when we
want to invoke this rule, we will just say ``by symmetry.''

The {\em crossed product\/} of $A$ by $\alpha$ is the \cstar completion, 
denoted $A \times_\alpha \Z$, of the {\em algebraic\/} direct sum (i.e., 
finite sums) $\bigoplus_n D_n$ equipped with ${}^*$-algebra structure
\begin{gather*}
(xy)_n = \sum_{k\in\Z} \alpha^k(\alpha^{-k}(x_k)y_{n-k}) \\
(x^*)_n = \alpha^n(x^*_{-n}).
\end{gather*}

\begin{defn}
$p_n$ denotes the identity element of the double-dual $D_n^{**}$, when the 
latter is canonically embedded as a weak* closed ideal of $A^{**}$. 
\end{defn}

Thus, the $p_n$ are central projections in $A^{**}$, and 
$D_n^{**}=A^{**}p_n$. By \eqref{mapideals} we have 
\begin{equation}
\label{mapprojections}
\alpha^n(p_kp_{-n}) = p_{n+k}p_n,
\end{equation}
where $\alpha$ has been canonically extended to a partial automorphism of 
$A^{**}$. 

A {\em covariant representation\/} of $(A,\alpha)$ on a Hilbert 
space $\H$ is a pair $(\pi,u)$, where $\pi$ is a 
nondegenerate representation of $A$ on $\H$ and $u \in \B(\H)$
satisfies
\begin{gather}
\label{rangeproj}
uu^* = \pi(p_1) \spacetext{and}
 u^*u = \pi(p_{-1}); \\
\label{covariance}
\ad u \circ \pi(a) = \pi \circ \alpha(a) \spacetext{for} a \in D_{-1}.
\end{gather}

Thus, $u$ is a partial isometry with range and domain projections 
$\pi(p_1)$ and $\pi(p_{-1})$, respectively. Similarly to what we did for 
partial automorphisms, we let $u^0=1$, and for $n>0$ we let $u^{-n}$ be the 
$n$th power of the adjoint $u^*$. 

\begin{lem}
If $(\pi, u)$ is a covariant representation, then $u^nu^{-n} = \pi(p_n)$ 
and $\ad u^n \circ \pi(a) = \pi \circ \alpha^n(a)$ for all $n \in \Z , a 
\in D_{-n}$. 
\end{lem}

\begin{proof}
This follows from the definitions and \eqref{mapprojections}.
\end{proof}

In particular, the $u^n$ are all partial isometries.

\begin{defn}
For a covariant representation $(\pi,u)$, we write
\[ C^*(\pi,u) = \overline{\sum_{n \in \Z} \pi(D_n)u^n}. \]
\end{defn}

A quick calculation shows that for $n,k \in \Z$, $a \in D_n$, and $b \in 
D_k$
\[ \pi(a)u^n\pi(b)u^k = \pi \circ \alpha^n(\alpha^{-n}(a)b)u^{n+k} \]
and
\[ (\pi(a)u^n)^* = \pi \circ \alpha^{-n}(a^*)u^{-n}, \]
so $C^*(\pi,u)$ is a \cstar algebra.
For every covariant representation $(\pi,u)$ of $(A,\alpha)$ there is a 
unique nondegenerate representation $\pi \times u$ of $A \times_\alpha 
\Z$ determined by
\[ \pi \times u(a) = \sum_n \pi(a_n)u^n
\spacetext{for } a \in \bigoplus_n D_n, \]
and conversely every nondegenerate representation of $A \times_\alpha \Z$ 
is associated in this way with a unique covariant representation.

$A$ is faithfully and nondegenerately embedded in $A \times_\alpha \Z$, 
hence in $(A \times_\alpha \Z)^{**}$; let $i\colon A \to (A \times_\alpha 
\Z)^{**}$ be this embedding. Then the universal representation of $A 
\times_\alpha \Z$ in $(A \times_\alpha \Z)^{**}$ is of the form $i \times 
m$ for a unique partial isometry $m \in (A \times_\alpha \Z)^{**}$. 
Moreover, $\alpha = \ad m$ and
\[ A \times_\alpha \Z = C^*(i,m) = \overline{\sum_n D_nm^n}. \]

\begin{defn}
We refer to $m$ as the {\em canonical partial isometry implementing 
$\alpha$ in $(A \times_\alpha \Z)^{**}$}.
\end{defn}

We now briefly review the elementary theory of coactions, specialized to 
$\Z$. A {\em coaction\/} of $\Z$ on a \cstar algebra $B$ is a nondegenerate 
injection $\delta\colon B \to B \otimes C^*(\Z)$ satisfying the {\em 
coaction identity\/} 
\[ (\delta \otimes \iota) \circ \delta
= (\iota \otimes \delta_\Z) \circ \delta, \]
where $\iota$ always denotes the identity automorphism, and
$\delta_\Z\colon C^*(\Z) \to C^*(\Z) \otimes C^*(\Z)$ is the
homomorphism determined 
by $\delta(n) = n \otimes n$ for $n \in \Z$.
Coactions of $\Z$ correspond bijectively to actions of $\hat\Z=\T$, and
if $\beta$ is the action of $\T$ associated to the coaction $\delta$ of
$\Z$, then the coaction identity for $\delta$ says exactly that
$\beta_s\beta_t=\beta_{st}$ for $s,t\in\T$. For example, if $\beta$ is
the action of $\T$ on $C(\T)$ given by
\[
\beta_s(f)(t)=f(st),
\]
then the associated coaction $\delta$ of $\Z$ is given on monomials by
\[
\delta(z^n)=z^n\otimes n\quad\text{for }z\in\Z.
\]
If $\delta$ is a coaction on $B$,
the {\em spectral subspaces\/} of $B$ are 
\[ B_n = \{ b \in B \mid \delta(b) = b \otimes n \}
\spacetext{for } n \in \Z. \]
The disjoint union of the $B_n$ forms a \cstar algebraic bundle over 
$\Z$:
\[ B_nB_k \subset B_{n+k} \spacetext{and } B^*_n = B_{-n}. \]

A {\em covariant representation\/} of $(B,\Z,\delta)$ is a pair 
$(\pi,\mu)$, where $\pi$ and $\mu$ are nondegenerate representations of $B$ 
and $c_0(\Z)$, respectively, satisfying a certain covariance condition. In 
this case it is convenient to recast the covariance as a relation between 
the spectral subspaces $B_n$ and the associated partition of unity 
\[ q_n = \mu(\chi_{\{n\}}). \]
By \cite[Lemma 2.2]{qui:disc}, the covariance condition becomes
\[ \pi(b)q_k = q_{n+k}\pi(b) \spacetext{for } n,k \in \Z, b \in B_n. \]
If every covariant representation factors through 
$(\pi,\mu)$, then $\overline{\pi(B)\mu(c_0(\Z))}$ is called the {\em 
crossed product\/} $B \times_\delta \Z$, and it is unique up to isomorphism. 

The coaction $\delta$ is called {\em inner\/} if there is a
nondegenerate homomorphism of $c_0(\Z)$ into $M(B)$ which is covariant
for the identity map of $B$, i.e., there is
a partition of unity $q_n$ in $M(B)$ such that $\lim_{n\to 
\pm\infty}q_n=0$ strictly and
\[ bq_k = q_{n+k}b \spacetext{for } n,k \in \Z, b \in B_n. \]
In this case by \cite[Theorem 6.9]{qui:twist} or \cite[Theorem 2.9]{lprs}
we have
\[ B \times_\delta \Z \cong B \otimes c_0(\Z). \]
In particular,
this holds for the {\em trivial\/} coaction $b \mapsto b \otimes 1$ of
$\Z$ on $B$ (take $q_n = 1$ for $n=0$ and $0$ otherwise).

A nondegenerate homomorphism $\rho$ of $B$ to $C$ is called {\em 
equivariant\/} for coactions $\delta$ and $\epsilon$ if
\[ \epsilon \circ \rho
= (\rho \otimes \iota) \circ \delta .\]
In this case, $\rho(B)$ is an $\epsilon$-{\em invariant\/} subalgebra of
$C$, i.e., $\epsilon(\rho(B)) \subset \rho(B) \otimes C^*(\Z)$.
If $I$ is an invariant ideal of $B$, $I \times_\delta \Z$ is an ideal of 
$B \times_\delta \Z$, and there is a natural coaction $\tilde\delta$ of 
$\Z$ on $B/I$ such that
\[ (B \times_\delta \Z) / (I \times_\delta \Z)
\cong (B/I) \times_{\tilde\delta} \Z. \]

If $\alpha$ is an actual automorphism of $A$ (so $n \mapsto \alpha^n$ is 
an action of $\Z$ on $A$), Takai-Takesaki duality says (in the language 
of coactions) that there is a coaction $\hat\alpha$ of $\Z$ on $A 
\times_\alpha \Z$ such that
\[ (A \times_\alpha \Z) \times_{\hat\alpha} \Z
\cong A \otimes \K(l^2(\Z)), \]
where $\K$ here stands for compact operators. While the construction of the 
dual coaction for partial automorphisms, indicated in the following 
proposition, is the same as for actions, the crossed product duality is 
largely destroyed, as we will see in \secref{duality}. 

\begin{prop} \cite{qui-rae:partial}
If $\alpha$ is a partial automorphism of $A$, then there is a unique 
coaction $\hat\alpha$ of $\Z$ on $A \times_\alpha \Z$ such that
\[ \hat\alpha(am^n) = am^n \otimes n
\spacetext{for } n \in \Z, a \in D_n. \]
The spectral subspaces are
\[ (A \times_\alpha \Z)_n = D_nm^n. \]
\end{prop}

\section{Invariant ideals}
\label{invariantideals}

\begin{defn}
A \cstar subalgebra $B$ of $A$ is called $\alpha$-{\em invariant\/} if
\begin{equation}
\label{invariantsubalgebra}
\alpha(B \cap D_{-1}) \subset B \spacetext{and }
\alpha\inv(B \cap D_1) \subset B.
\end{equation}
\end{defn}

\begin{lem}
If \textup(and only if\textup) $B$ is $\alpha$-invariant, then
\begin{equation}
\label{invariantequal}
\alpha^n(B \cap D_{-n}) = B \cap D_n \spacetext{for } n \in \Z.
\end{equation}
\end{lem}

\begin{proof}
\eqref{invariantequal} is trivial for $n=0$. Since $\alpha$ and 
$\alpha\inv$ are inverses, \eqref{invariantequal} for $n=\pm 1$ follows 
from \eqref{invariantsubalgebra}. Inductively, assume $k>1$ and
\eqref{invariantequal} holds for $|n| < k$. Then
\begin{align*}
\begin{split}
\alpha^k(B \cap D_{-k})
&= \alpha(\alpha^{k-1}(B \cap D_{-k}))
\subset \alpha(B \cap D_{k-1} \cap D_{-1}) \\
&\subset B \cap D_1 \cap D_k
= B \cap D_k.
\end{split}
\end{align*}
Symmetrically,
\[ \alpha^{-k}(B \cap D_k) \subset B \cap D_{-k}. \]
Since $\alpha^k$ and $\alpha^{-k}$ are inverses, we must have 
\eqref{invariantequal} for $n=\pm k$.
\end{proof}

Thus, if $B$ is $\alpha$-invariant, then $\alpha$ restricts to a
partial automorphism, which we also denote by $\alpha$, of $B$, with
ideals $B \cap D_n$.

McClanahan \cite[Proposition 5.1 and Corollary 5.2]{mcc} proves most of the 
following result for partial actions of arbitrary discrete groups; the 
basis for the techniques of Propositions \ref{idealcrossedproduct} and 
\ref{quotientcrossedproduct} can actually be found in \cite[Lemma 
1]{gre:smooth}. Since our notation is different, and since we include the 
dual coaction, we give the proof for the reader's convenience. 

\begin{prop}
\label{idealcrossedproduct}
Let $I$ be an $\alpha$-invariant ideal of $A$, let $q$ be the identity 
element of $I^{**}$ in $A^{**}$, and let $i\colon I \hookrightarrow A$ 
be the inclusion map. Then $(i,qm)$ is a covariant representation of 
$(I,\alpha)$, and $i \times qm$ is an isomorphism of $I \times_\alpha \Z$ 
onto the ideal $\overline{\sum_n ID_nm^n}$ of $A \times_\alpha \Z$. 
Moreover, this isomorphism is equivariant with respect to the dual 
coactions.
\end{prop}

\begin{proof}
We first show that $m$ commutes with $q$ in $(A \times_\alpha \Z)^{**}$:
\begin{align*}
\begin{split}
qm &= qp_1m = \alpha(qp_{-1})m = mqp_{-1}m^*m \\
&= mqp_{-1}^2 = mp_{-1}q = mq.
\end{split}
\end{align*}
Thus,
\[ (qm)^n(qm)^{-n} = qm^nm^{-n} = qp_n \spacetext{for } n \in \Z, \]
and, since  $q$ is the identity of $I^{**}$,
\[ \ad (qm) \circ i(a) = \ad m \circ i(a) = \alpha(a)
\spacetext{for } a \in ID_{-1}. \]
We have
\[ (i \times qm)(I \times_\alpha \Z)
= \overline{\sum_n ID_n(qm)^n}
= \overline{\sum_n ID_nm^n}, \]
which is obviously an ideal of $A \times_\alpha \Z$. The equivariance is 
now obvious.

It remains to show $i \times qm$ is injective, and this is accomplished by 
showing that every covariant representation $(\pi,u)$ of $(I,\alpha)$ 
factors through $(i,qm)$. Let $\bar\pi$ be the unique representation of $A$ 
extending $\pi$. We verify that $(\bar\pi,u)$ is a covariant representation 
of $(A,\alpha)$. First, if $\{e_i\}$ is an approximate identity of $I$, 
then (taking weak operator limits) 
\[ \bar\pi(p_n) = \lim \pi(e_ip_n) = \pi(qp_n) = u^nu^{-n}. \]

Now change $\{e_i\}$ to an approximate identity of $ID_{-1}$. Then 
$\{\alpha(e_i)\}$ is an approximate identity of $ID_1$, so for $a \in 
ID_{-1}$ we have
\begin{align*}
\begin{split}
\ad u \circ \bar\pi(a)
&= \lim \ad u \circ \pi(e_ia)
= \lim \pi \circ \alpha(e_ia) \\
&= \lim \pi(\alpha(e_i)\alpha(a))
= \bar\pi \circ \alpha(a). \hfill\qed
\end{split}
\end{align*}
\noqed
\end{proof}

Let $I$ be an $\alpha$-invariant ideal of $A$, and let $\zeta \colon A \to 
A/I$ be the quotient map. Then $\tilde\alpha \circ \zeta = \zeta \circ 
\alpha$ determines a partial automorphism $\tilde\alpha$ of $A/I$, with 
ideals $\zeta(D_n)$. Moreover, the identity of $\zeta(D_n)^{**}$ is 
$\zeta(p_n)$ (more precisely, $\zeta^{**}(p_n)$).

\begin{prop}
\label{quotientcrossedproduct}
Let $I$, $\zeta$, and $\tilde\alpha$ be as above, and let $\tilde m$ be the 
canonical partial isometry implementing $\tilde\alpha$ in $(A/I 
\times_{\tilde\alpha} \Z)^{**}$. Then $(\zeta,\tilde m)$ is a covariant 
representation of $(A,\alpha)$, and $\zeta \times \tilde m$ is a surjection 
of $A \times_\alpha \Z$ onto $A/I \times_{\tilde\alpha} \Z$ with kernel $I 
\times_\alpha \Z$. Moreover, this surjection is equivariant with respect to 
the dual coactions. 
\end{prop}

\begin{proof}
We have
\[ \tilde m^n \tilde m^{-n} = \zeta(p_n) \spacetext{for } n \in \Z, \]
and
\[ \ad \tilde m \circ \zeta = \tilde\alpha \circ \zeta
= \zeta \circ \alpha. \]
Clearly $\zeta \times \tilde m$ is a surjection of $A \times_\alpha \Z$ 
onto $A/I \times_{\tilde\alpha} \Z$.

Since $\zeta \times \tilde m$ vanishes on $ID_nm^n$ for every $n$, 
$\ker(\zeta \times \tilde m) \supset I \times_\alpha \Z$. For the opposite 
containment, let $(\pi,u)$ be a covariant representation of $(A,\alpha)$ 
with $\ker(\pi \times u) = I \times_\alpha \Z$. Then $\ker\pi \supset I$ 
since $I\times_\alpha\Z\supset I$ and $(\pi\times u) | I = \pi | I$. So, 
there is a representation $\tilde\pi$ of $A/I$ such that $\pi = \tilde\pi 
\circ \zeta$. Then $(\tilde\pi,u)$ is a covariant representation of 
$(A/I,\tilde\alpha)$, and
\[ \pi \times u = (\tilde\pi \times u) \circ (\zeta \times \tilde m). \]
Hence
\[ \ker(\zeta \times \tilde m)
\subset \ker(\pi \times u) = I \times_\alpha \Z. \]

For the equivariance, if $a \in D_n$ we have
\begin{align*}
\begin{split}
\hat{\tilde\alpha} \circ (\zeta \times \tilde m)(am^n)
&= \hat{\tilde\alpha}(\zeta(a) \tilde m^n)
= \zeta(a) \tilde m^n \otimes n \\
&= ((\zeta \times \tilde m) \otimes \iota)(am^n \otimes n) \\
&= ((\zeta \times \tilde m) \otimes \iota) \circ \hat\alpha(am^n). \hfill\qed
\end{split}
\end{align*}
\noqed
\end{proof}

We will need the following elementary result on tensor products of 
partial automorphisms.

\begin{prop}
Let $\alpha \in \paut A$, and let $B$ be a \cstar algebra. Then $\iota 
\otimes \alpha \in \paut B \otimes A$, with domain $B \otimes D_{-1}$, and 
\[ (B \otimes A) \times_{\iota \otimes \alpha} \Z \cong B \otimes (A 
\times_\alpha \Z), \] where the minimal tensor product is used throughout. 
\end{prop}

\begin{proof}
Crossed products by partial automorphisms of $\Z$ are automatically reduced 
\cite[Theorem 5.2]{exe:partial}, \cite[Proposition 4.2]{mcc}. If $B$ and 
$A$ are faithfully and nondegenerately represented on Hilbert spaces $\K$ 
and $\H$, respectively, then $B \otimes A$ is so represented on $\K \otimes 
\H$. So, $(B \otimes A) \times_{\iota \otimes \alpha,r} \Z$ and $B \otimes 
(A \times_{\alpha,r} \Z)$ are both represented on $\K \otimes \H \otimes 
l^2(G)$, and a mildly careful examination of these representations yields 
the fruit that these \cstar algebras are in fact equal. 
\end{proof}

\section{duality}
\label{duality}

In this section we begin to examine the extent to which Takai-Takesaki 
crossed product duality fails for partial automorphisms. 

\begin{defn}
Let
\[ D_\infty = \bigcap_{n>0} D_n \spacetext{and}
 D_{-\infty} = \bigcap_{n<0} D_n .\]
\end{defn}

\begin{lem}
$D_\infty$ and $D_{-\infty}$ are $\alpha$-invariant.
\end{lem}

\begin{proof}
Let $a \in D_\infty D_{-1}$. Then $a \in D_n D_{-1}$ for all $n>0$, so 
$\alpha(a) \in D_{n+1}$ for all $n>0$, hence $\alpha(a) \in D_\infty$. If 
$a \in D_\infty$, then $a \in D_{n+1}$ for all $n>0$, so $\alpha\inv(a) 
\in D_n$ for all $n>0$, so $\alpha\inv(a) \in D_\infty$. The invariance of 
$D_{-\infty}$ follows by symmetry.
\end{proof}

\begin{prop}
$\alpha$ restricts to an automorphism of $D_\infty D_{-\infty}$.
\end{prop}

\begin{proof}
Since both $D_\infty$ and $D_{-\infty}$ are $\alpha$-invariant, so is 
$D_\infty D_{-\infty}$. Since $D_\infty D_{-\infty} \subset D_1D_{-1}$, we 
are done. 
\end{proof}

The next result shows that $D_\infty D_{-\infty}$ is a kind of 
``automorphic core'' of $\alpha$.

\begin{lem}
$D_\infty D_{-\infty} = \{0\}$ if and only if there is no nonzero \cstar 
subalgebra $B$ of $A$ such that $\alpha | B$ is an automorphism.
\end{lem}

\begin{proof}
If such a $B$ exists, we have $B \subset D_n$ for all $n$, so
\[ D_\infty D_{-\infty} \supset B \ne \{0\}. \]

The converse is trivial since $\alpha | D_\infty D_{-\infty}$ is an 
automorphism.
\end{proof}

\begin{defn}
We call $\alpha$ {\em completely nonautomorphic\/} if $D_\infty 
D_{-\infty} = \{0\}$. 
\end{defn}

\begin{prop}
The quotient partial automorphism on $A/(D_\infty D_{-\infty})$ is 
completely nonautomorphic.
\end{prop}

\begin{proof}
Let $\zeta\colon A \to A/(D_\infty D_{-\infty})$ be the quotient map. 
Recall that the ideals of the quotient partial automorphism are 
$\zeta(D_n)$. Since $D_\infty D_{-\infty} \subset D_n$ for all $n$, we have
\begin{align*}
\begin{split}
\bigcap_n \zeta(D_n)
&= \bigcap_n D_n/(D_\infty D_{-\infty})
= \bigl( \bigcap_n D_n \bigr) / (D_\infty D_{-\infty}) \\
&= (D_\infty D_{-\infty}) / (D_\infty D_{-\infty})
= \{0\}. \hfill\qed
\end{split}
\end{align*}
\noqed
\end{proof}

The following result shows that crossed product duality tends to be 
maximally false for completely nonautomorphic partial automorphisms. 
Let $p_{\infty}$ (respectively, $p_{-\infty}$) denote the identity 
of $D^{**}_{\infty}$ (respectively, $D^{**}_{-\infty}$) in $A^{**}$.

\begin{thm}
\label{nonautomorphicinner}
If $\alpha$ is completely nonautomorphic and 
$\lim_{n\to\pm\infty}p_{n}=p_{\pm\infty}$ strictly in $M(A)$, then the 
dual coaction $\hat\alpha$ of $\Z$ on $A \times_\alpha \Z$ is inner, 
so \[ (A \times_\alpha \Z) \times_{\hat\alpha} \Z \cong (A 
\times_\alpha \Z) \otimes c_0(\Z).  \]
\end{thm}

\begin{proof}
We must produce a partition of unity $\{ q_n \mid n \in \Z \}$ in $M(A 
\times_\alpha \Z)$ such that $\lim_{n\to\pm\infty}q_{n}=0$ strictly 
and
\begin{equation}
\label{inner}
am^nq_k = q_{n+k}am^n \spacetext{for } n,k \in \Z, a \in D_n.
\end{equation}
Define
\[ q_n = \begin{cases}
p_\infty(p_{n+1} - p_n) & \text{if } n<0, \\
p_n - p_{n+1} & \text{if } n\ge 0.
\end{cases} \]

The $q_n$ for $n\ge 0$ form a partition of $1 - p_\infty$, and
the $q_n$ for $n<0$ form a 
partition of $p_\infty$, since
\[ 0 = p_\infty p_{-\infty}
= \text{weak*-} \lim_{n \to -\infty} p_\infty p_n. \]
Further, the hypotheses imply $\lim_{n\to\pm\infty}q_{n}=0$ 
strictly in $M(A)$, hence in $M(A\times_{\alpha}\Z)$.

By induction and taking adjoints, and since $m^0=1$ and the $q_n$ are 
in the center of $A^{**}$, \eqref{inner} will follow from
\begin{equation}
\label{mqinner}
mq_k = q_{k+1}m \spacetext{for } k \in \Z.
\end{equation}
We first show
\begin{equation}
\label{mpinner}
mp_k = p_{k+1}m \spacetext{for } k \in \Z.
\end{equation}
For $k\ge 0$
\begin{align*}
\begin{split}
mp_k
&= mp_{-1}p_k
= mp_kp_{-1}
= mm^km^{-k}m^{-1}m \\
&= m^{k+1}m^{-k-1}m
= p_{k+1}m,
\end{split}
\end{align*}
while for $k<0$
\begin{align*}
\begin{split}
mp_k
&= mm^km^{-k}
= mm^{-1}m^{k+1}m^{-k-1}m
= p_1p_{k+1}m \\
&= p_{k+1}p_1m
= p_{k+1}m,
\end{split}
\end{align*}
showing \eqref{mpinner}. Since $m$ commutes with $p_\infty$, 
\eqref{mqinner} follows for all $k\ne -1$. For the remaining case,
\begin{align*}
\begin{split}
mq_{-1}
&= mp_\infty(p_0 - p_{-1})
= p_\infty m(p_0 - p_{-1}) \\
&= p_\infty (mp_0 - mp_{-1})
= p_\infty (m - m)
= 0 \\
&= p_0m - p_1m
= (p_0 - p_1)m
= q_0m. \hfill\qed
\end{split}
\end{align*}
\noqed
\end{proof}

\begin{cor}
\label{nonautomorphicinnercor}
Let $\alpha$ be completely nonautomorphic, and assume the projections 
$p_{\pm 1}$ and $p_{\pm\infty}$ are in $M(A)$.  Then \textup(assuming 
$A\ne\{0\}$\textup) $A$ contains a nonzero $\alpha$-invariant ideal 
$I$ such that $\hat\alpha$ is inner on $I\times_{\alpha}\Z$.
\end{cor}

\begin{proof}
First note that $m\in M(A\times_{\alpha}\Z)$ by \cite[Proposition 
2.13]{qui-rae:partial}, so for all $n\in Z$ we have $m^{n}\in 
M(A\times_{\alpha}\Z)$, hence $p_{n}\in M(A)$.
Let
\begin{align*}
I_{1}&=\overline{\sum_{n<0}p_{\infty}(p_{n+1}-p_{n})A}\\
I_{2}&=\overline{\sum_{n>0}p_{-\infty}(p_{n-1}-p_{n})A}\\
I_{3}&=\overline{\sum_{n,k>0}(1-p_{\infty}-p_{-\infty})
(p_{1-n}-p_{-n})(p_{k-1}-p_{k})A}
\end{align*}
Then each $I_{j}$ is an $\alpha$-invariant ideal to which the 
hypotheses of \thmref{nonautomorphicinner} apply, and 
$A=\bigoplus_{1}^{3}I_{j}$.
\end{proof}

\section{The distribution of ideals}
\label{distribution}

The ideals $D_n$ decrease in both directions from $n=0$. In each direction, 
the behavior of the partial action is dramatically influenced by whether 
the ideals are eventually constant or strictly decreasing forever, and upon 
whether the intersection is $\{0\}$. 

\begin{lem}
\label{nilpotent}
For $n>0$, the following are equivalent:
\begin{enumerate}
\item $D_n = \{0\}$\textup;
\item $D_{-n} = \{0\}$\textup;
\item $D_{n-1}D_{-1} = \{0\}$.
\item $D_{1-n}D_1 = \{0\}$\textup;
\end{enumerate}
\end{lem}

\begin{proof}
Since $\alpha$ is injective, this follows from the following relations:
\[ D_n = \alpha^n(D_{-n}) = \alpha(D_{n-1}D_{-1})
= \alpha^{n-1}(D_{1-n}D_1). \hfill\qed \]
\noqed
\end{proof}

\begin{defn}
A partial automorphism such that $D_n = \{0\}$ for some $n>0$ will be 
called {\em nilpotent}. 
\end{defn}

The next result concerns the most trivial nilpotent partial automorphisms. 

\begin{prop}
\label{trivial}
If $D_1 = \{0\}$, then $A \times_\alpha \Z = A$ and the dual coaction 
$\hat\alpha$ is the trivial coaction $a \mapsto a \otimes 1$, so
\[ (A \times_\alpha \Z) \times_{\hat\alpha} \Z
\cong A \otimes c_0(\Z). \]
\end{prop}

\begin{proof}
This follows immediately from the definitions, since $\bigoplus_n D_n = A$ 
and $m = 0$.
\end{proof}

Perhaps the simplest nontrivial nilpotent partial automorphisms are given 
by the following example. 

\begin{ex}
Define $\sigma_n \in \paut \C^n$ by
\[ \sigma_n(z_1,\dots ,z_{n-1},0) = (0,z_1,\dots ,z_{n-1}). \]
We will refer to this as the {\em shift on $\C^n$}. Exel 
\cite{exe:partial} shows that
\[ C^n \times_{\sigma_n} \Z \cong M_n, \]
the algebra of $n \times n$ matrices.
\end{ex}

The next lemma shows that all nilpotent partial automorphisms contain 
shifts.

\begin{lem}
Let $n>1$, and suppose $D_{n-1} \ne \{0\}$, $D_n = \{0\}$, and
\[ A = D_{1-n} + D_{2-n}D_1 + D_{3-n}D_2 + \dots + D_{n-1}. \]
Then $(A,\alpha)$ is isomorphic to $(D_{n-1} \otimes \C^n,\iota \otimes 
\sigma_n)$, so
\[ A \times_\alpha \Z \cong D_{n-1} \otimes M_n. \]
\end{lem}

\begin{proof}
It suffices to prove $(A,\alpha)$ is isomorphic to $(D_{1-n} \otimes 
\C^n,\iota \otimes \sigma_n)$, since composing with $\alpha^{n-1} \otimes 
\iota$ will then give an isomorphism with $(D_{n-1} \otimes \C^n,\iota 
\otimes \sigma_n)$. The hypothesis implies 
\[ D_{-1} = D_{1-n} + D_{2-n}D_1 + D_{3-n}D_2 + \dots + D_{-1}D_{n-2} \]
and that the ideals $D_{1-n}$, $D_{2-n}D_1$, $D_{3-n}D_2$, \dots, $D_{n-1}$ 
have pairwise zero intersection.
Define $\theta\colon D_{1-n} \otimes \C^n \to A$ by 
\[ \theta(a \otimes (z_1,\dots ,z_n)) = \sum_1^n z_i\alpha^{i-1}(a). \]
$\theta$ is clearly an isomorphism, and
\[ \theta(D_{1-n} \otimes (\C^{n-1} \times \{0\}))
= \sum_1^{n-1} \alpha^{i-1}(D_{1-n})
= \sum_1^{n-1} D_{i-n}D_{i-1} = D_{-1}. \]
We have
\begin{align*}
\theta \circ (\iota \otimes \sigma_n)(a \otimes (z_1,\dots ,z_{n-1},0)
&= \theta(a \otimes (0,z_1,\dots ,z_{n-1})) \\
&= \sum_1^{n-1} z_i\alpha^i(a) \\
&= \alpha(\sum_1^{n-1} z_i\alpha^{i-1}(a)) \\
&= \alpha \circ \theta(a \otimes (z_1,\dots ,z_{n-1},0)). \hfill\qed
\end{align*}
\noqed
\end{proof}

The following theorem shows that $A$ can have many subquotients (in 
fact, is sometimes an inverse limit of such) on which $\alpha$ looks 
like a nilpotent shift.

\begin{thm}
\label{subquotient}
For each $n>1$ the ideal
\[ I_n = D_{1-n} + D_{2-n}D_1 + D_{3-n}D_2 + \dots + D_{n-1} \]
of $A$ is $\alpha$-invariant. Moreover, $I_n \supset I_{n+1}$, and if 
$\beta_n$ is the quotient partial automorphism of $I_n/I_{n+1}$, then
\[ (I_n/I_{n+1},\beta_n)
\cong (D_{n-1} / (D_n + D_{n-1}D_{-1}) \otimes \C^n,
\iota \otimes \sigma_n). \]
Consequently,
\[ I_n/I_{n+1} \times_{\beta_n} \Z
\cong D_{n-1} / (D_n + D_{n-1}D_{-1}) \otimes M_n. \]
\end{thm}

\begin{proof}
We have
\begin{align*}
\begin{split}
\alpha(I_nD_{-1})
&= \alpha(D_{1-n} + D_{2-n}D_1 + \dots + D_{-1}D_{n-2}) \\
&= D_{2-n}D_1 + \dots + D_{n-1} \subset I_n.
\end{split}
\end{align*}
The containment $I_n \supset I_{n+1}$ follows from $D_n \supset D_{n+1}$. 
Noting that the ideals of $\beta_n$ are
\[ D_k(\beta_n) = I_nD_k / I_{n+1}D_k, \]
we see that $(I_n/I_{n+1},\beta_n)$ satisfies the hypotheses of the above 
lemma. We finish by observing that
\[ I_{n+1}D_{n-1} = D_n + D_{n-1}D_{-1}, \]
so
\[ D_{n-1}(\beta_n) = D_{n-1}/(D_n + D_{n-1}D_{-1}). \hfill\qed \]
\noqed
\end{proof}

A version of the above result holds even for $n=1$, although it 
(like the empty set) is best dealt with separately:

\begin{thm}
\label{subquotientspecial}
Let $I = D_{-1} + D_1$. Then
\[ (A \times_\alpha \Z) / (I \times_\alpha \Z) \cong A/I \]
and
\[ \bigl( (A \times_\alpha \Z) \times_{\hat\alpha} \Z \bigr) /
\bigl( (I \times_\alpha \Z) \times_{\hat\alpha} \Z \bigr)
\cong A/I \otimes c_0(\Z). \]
\end{thm}

\begin{proof}
We have
\[ (A \times_\alpha \Z) / (I \times_\alpha \Z)
\cong A/I \times_\beta \Z  \]
and
\[ \bigl( (A \times_\alpha \Z) \times_{\hat\alpha} \Z \bigr) /
\bigl( (I \times_\alpha \Z) \times_{\hat\alpha} \Z \bigr)
\cong (A/I \times_{\beta} \Z) \times_{\hat\beta} \Z, \]
where $\beta$ is the quotient partial automorphism.
But $\beta$ has domain $\{0\}$, so by \propref{trivial}
\[ A/I \times_\beta \Z = A/I \]
and
\[ (A/I \times_{\beta} \Z) \times_{\hat\beta} \Z
\cong A/I \otimes c_0(\Z). \]
\end{proof}

We use the above results to show that crossed product duality
fails in general for partial automorphisms.  To be precise, crossed 
product duality demands a {\em canonical\/} isomorphism of $(A 
\times_\alpha \Z) \times_{\hat\alpha} \Z$ with $A \otimes 
\K(l^2(\Z))$---an ``accidental'' isomorphism is irrelevant.  Without 
putting too fine a point on it, let us agree that ``canonical'' 
implies at least that if $I\supset J$ are $\alpha$-invariant ideals of 
$A$, then the isomorphism carries $(I \times_\alpha \Z) 
\times_{\hat\alpha} \Z$ onto $I \otimes \K(l^2(\Z))$, and similarly 
for $J$.  Taking quotients, we get an isomorphism of $(I/J 
\times_{\tilde\alpha} \Z) \times_{\hat{\tilde\alpha}} \Z$ with $I/J 
\otimes \K(l^2(\Z))$.  But with $I=I_{n}$ and $J=I_{n+1}$, the above 
results would then imply
\[ I/J \otimes \K(l^2(\Z)) \cong I/J\otimes M_{n} \otimes c_0(\Z), \]
which is false except in very special examples.
This failure of crossed product 
duality has an advantage over \corref{nonautomorphicinnercor}, since 
it does not require any projections to be multipliers.  However, this 
does not handle all cases: the following example shows that even when 
$D_{\infty}=D_{-\infty}=\{0\}$, we can have $I_{n}=A$ for all $n>0$.

\begin{ex}
This example was invented by N\'andor Sieben. Here $A$ will be 
$C_{0}(\R)$ and $\alpha$ will be translation by $\pi$:
\[\alpha(f)(t)=f(t-\pi),\]
with domain
\[D_{-1}=\left\{f\in A\mid f(t)=0\text{ for }t\in S\right\},\]
where
\[S=\{0\}\cup\{\pm\sum_{1}^{n}\frac{1}{k}\mid n\in\N\}.\]
The reader can check that
\[\bigcup_{n>0}(S+n\pi)\quad\text{and}\quad
\bigcup_{n<0}(S+n\pi)\]
are both dense in $\R$, so
\[D_{\infty}=D_{-\infty}=\{0\},\]
and moreover for each $n>0$
\[\bigcup^{n}_{k=1}(S+k\pi)\cap\bigcup^{n-1}_{k=0}(S-k\pi)=\emptyset,\]
so
\[I_{n+1}\supset D_{n}+D_{-n}=A.\]
\end{ex}

We saw in the preceding section that $D_\infty$ and $D_{-\infty}$ are 
$\alpha$-invariant. What can we say about the restricted and quotient 
partial automorphisms?

\begin{defn}
$\alpha$ is a {\em forward shift\/} (respectively, a {\em backward 
shift\/}) if $D_n=A$ for all $n\le 0$ and 
$D_\infty = \{0\}$ (respectively, $D_n=A$ for all $n\ge 0$ and 
$D_{-\infty} = \{0\}$).
\end{defn}

Of course, $\alpha$ is a forward shift if and only if $\alpha\inv$ is 
a backward shift.  The simplest nonnilpotent forward shift is
\[ \sigma\colon (x_1,x_2,\dotsc) \mapsto (0,x_1,x_2,\dotsc) \]
on $c_0$, whose crossed product is the compact operators.  If we 
adjoin an identity to $c_{0}$, the crossed product becomes the 
Toeplitz algebra generated by a nonunitary isometry 
\cite{exe:partial}.  The following result shows that nonnilpotent 
forward shifts tend to look like the preceding example on a certain 
ideal.

\begin{prop}
\label{shiftiso}
Let $\alpha$ be a forward shift, and assume $p_{1}\in M(A)$. Then
\[I=\overline{\sum_{n>0}(p_{n-1}-p_{n})A}\]
is and $\alpha$-invariant ideal, and
\[(I,\alpha|I)\cong ((1-p_{1})A\otimes c_{0},\iota\otimes\sigma),\]
where $\sigma$ is the forward shift on $c_{0}$. A similar result 
holds for backward shifts.
\end{prop}

\begin{proof}
The hypotheses imply
\[I=\bigoplus_{n\ge 0}\alpha^{n}((1-p_{1})A),\]
and the proposition follows easily.
\end{proof}

\begin{prop}
If $\alpha$ is completely nonautomorphic, then $\alpha | D_{-\infty}$ 
is a forward shift and $\alpha|D_\infty$ is a backward shift.
\end{prop}

\begin{proof}
The ideals for $\alpha | D_{-\infty}$ are $D_nD_{-\infty}$, which 
coincide with $D_{-\infty}$ for $n\le 0$, and we have
\[ \bigcap_{n>0} D_nD_{-\infty} = D_\infty D_{-\infty} = \{0\}. \]
The other statement follows by symmetry.
\end{proof}

Note that it is possible for $D_\infty D_{-\infty} = \{0\}$ while 
neither $D_\infty$ nor $D_{-\infty}$ is $\{0\}$, e.g., the direct sum 
of a forward shift and a backward shift.  In fact, this is almost 
typical, as we will discuss after the next result.

\begin{prop}
If $\tilde\alpha$ is the quotient partial automorphism on 
$A/D_\infty$ \textup(respectively, $A/D_{-\infty}$, 
$A/(D_{\infty}+D_{-\infty})$\textup), 
then $D_\infty(\tilde\alpha) = \{0\}$ \textup(respectively, 
$D_{-\infty}(\tilde\alpha) = \{0\}$,
$D_\infty(\tilde\alpha) = D_{-\infty}(\tilde\alpha) =\{0\}$\textup).
\end{prop}

\begin{proof}
The first part follows from
\[ (D_n + D_\infty)/D_\infty
= D_n/D_\infty \spacetext{for} n>0, \]
and the other parts are shown similarly.
\end{proof}

In particular, a completely nonautomorphic partial automorphism falls 
naturally into three pieces: a forward shift on $D_{-\infty}$, a 
backward shift on $D_{\infty}$, and a quotient partial automorphism on 
$A/(D_{\infty}+D_{-\infty})$ satisfying 
$D_{\infty}=D_{-\infty}=\{0\}$.  The simplest nonnilpotent 
illustration of the latter phenomenon is 
$\bigoplus_{n>0}(\C^{n},\sigma^{n})$.  This can also be visualized as 
the partial automorphism on $c_{0}(\N^{2})$ with domain $\{x\mid 
x_{n,1}=0\text{ for all }n\in\N\}$ and which takes such an $x$ to $y$, 
where
\[y_{n,k}=
\begin{cases}
x_{n-1,k+1}&\text{if $n>1$,}\\
0&\text{if $n=1$.}
\end{cases}\]
Note that in this case
\[I_{n}=\{x_{k,l}\mid k+l\ge n+1\}.\]
The following result shows that partial automorphisms with 
$D_{\infty}=D_{-\infty}=\{0\}$ tend to look like the preceding 
example on a certain ideal.

\begin{prop}
Let $\alpha$ satisfy $D_{\infty}=D_{-\infty}=\{0\}$, and assume 
$p_{1},p_{-1}\in M(A)$. Then
\[I=\overline{\sum_{n,k>0}(p_{1-n}-p_{-n})(p_{k-1}-p_{k})A}\]
is an $\alpha$-invariant ideal, and
\[(I,\alpha|I)\cong
\bigoplus_{n>0}((p_{1-n}-p_{-n})(1-p_{1})A\otimes\C^n,
\iota\otimes\sigma_{n}).\]
\end{prop}

\begin{proof}
The proof is almost as easy as in \propref{shiftiso}, noting that
\[\alpha((p_{-k}-p_{-k-1})(p_{l-1}-p_{l}))
=(p_{1-k}-p_{-k})(p_{l}-p_{l+1})\quad\text{for}\quad k,l>0\]
and the hypotheses imply
\[I=\sum^{\infty}_{n=1}\sum^{n-1}_{k=0}
\alpha^{k}((p_{1-n}-p_{-n})(1-p_{1})A).\hfill\qed\]
\noqed
\end{proof}

Note that the summands in the above proposition are the subquotients 
of Theorems \ref{subquotient} and \ref{subquotientspecial}.

The nilpotent partial automorphisms studied in the preceding section 
gave trivial examples of the $D_n$ being eventually constant.  We 
finish by examining the general case.

\begin{lem}
If $n\le 0$ and $D_n = D_{n-1}$, then
\begin{equation}
\label{constant}
D_n = D_k \spacetext{for all } k \le n,
\end{equation}
and similarly for $n\ge 0$.
\end{lem}

\begin{proof}
The hypothesis implies
\[ D_{n-1} = \alpha\inv(D_nD_1) = \alpha\inv(D_{n-1}D_1) = D_{n-2}, \]
giving \eqref{constant} by induction. The other part follows by symmetry.
\end{proof}

Thus, if $n\le 0$ and $D_n = D_{n-1}$, then $D_n = D_{-\infty}$ is 
$\alpha$-invariant. Curiously, a partial converse holds:

\begin{prop}
If $n<0$ and $\alpha(D_n) \subset D_n$, then $D_n = D_{n-1}$.
\end{prop}

\begin{proof}
We have
\begin{align*}
\begin{split}
D_nD_{n-1}
&= \alpha\inv(D_{n+1}D_nD_1)
= \alpha\inv(\alpha(D_n)D_n) \\
&= \alpha\inv \circ \alpha(D_n)
= D_n,
\end{split}
\end{align*}
so $D_n \subset D_{n-1}$, whence $D_n = D_{n-1}$.
\end{proof}

What can we say about $\alpha | D_n$? It depends on 
the $D_k$ for $k\ge 0$:

\begin{lem}
If
\begin{align*}
n&=\max\{j\le 0\mid D_j=D_{j-1}\}\quad\text{and}\\
k&=\min\{j\ge 0\mid D_j=D_{j+1}\}
\end{align*}
are both finite,
then $n=-k$ 
and $\alpha$ restricts to an automorphism of $D_n$. 
\end{lem}

\begin{proof}
Assuming without loss of generality that $n\ge -k$, it suffices to show 
$\alpha(D_n)=D_n$:
\begin{align*}
\begin{split}
\alpha(D_n)
&= \alpha(D_{n-1})
= D_nD_1
= D_{-k}D_1 \\
&= \alpha^{-k}(D_kD_{k+1})
= \alpha^{-k}(D_k)
= D_{-k}
= D_n. \hfill\qed
\end{split}
\end{align*}
\noqed
\end{proof}

Thus, if the $D_n$ are eventually constant in both directions, then this 
behavior starts at the same place forward and backward, generalizing 
\lemref{nilpotent}. Moreover, when $D_n = D_{n-1} = D_k = D_{k+1}$ for $n 
= -k \le 0$, then $D_n$ is the automorphic core $D_\infty D_{-\infty}$. 
On the other hand, if $n\le 0$ and $D_n = D_{n-1}$, but $D_k \ne 
D_{k+1}$ for all $k>0$, then in any event $D_\infty \subset D_n$ and the 
quotient partial automorphism on $A/D_n$ is nilpotent. As usual, similar 
statements hold for $n\ge 0$.

\providecommand{\bysame}{\leavevmode\hbox to3em{\hrulefill}\thinspace}

\end{document}